\documentclass[12pt,a4paper]{article}
\usepackage[cp1251]{inputenc} % cp866 - DOS, cp1251 - Windows
\usepackage{amsfonts}

\newcommand{\di}{\displaystyle}
\newcommand{\B}{$\hfill\Box$}

\newcommand{\ga}{\gamma}
\newcommand{\de}{\delta}
\newcommand{\la}{\lambda}
\newcommand{\om}{\omega}
\newcommand{\ee}{\varepsilon}
\newcommand{\vv}{\varphi}

\begin{document}

\begin{center}
{\large\bf
Inverse Problems for Systems of Variable Order Differential Equations with 
Singularities on Spatial Networks.}\\[0.2cm]
{\bf V.\,Yurko} \\[0.2cm]
\end{center}

\thispagestyle{empty}

{\bf Abstract.} Variable order differential equations with non-integrable 
singularities are considered on spatial networks. Properties of the spectrum 
are established, and the solution of the inverse spectral problem is obtained.



{\bf 1. } 
Consider a compact star-type graph $T$ in ${\bf R^\om}$ with the set 
of vertices $V=\{v_0,\ldots, v_p\}$ and the set of edges ${\cal E}=
\{e_1,\ldots, e_p\},$ where $v_1,\ldots, v_{p}$ are the boundary vertices, 
$v_0$ is the internal vertex, and $e_j=[v_{j},v_0],$ $e_1\cap\ldots\cap 
e_p=\{v_0\}$. Let $l_j$ be the length of the edge $e_j$. For each edge 
$e_j\in {\cal E}$ we introduce the parameter $x_j\in [0,l_j]$ such that 
$x_j=0$ corresponds to the boundary vertices $v_1,\ldots, v_{p}$.
A function $Y$ on $T$ may be represented as $Y=\{y_j\}_{j=\overline{1,p}}$,
where the function $y_j(x_j)$ is defined on the edge $e_j$.

Let $n_j$, $j=\overline{1,p},$ be positive integers such that
$n_1\ge n_2\ge\ldots\ge n_p\ge 2.$ Consider the differential
equations on $T$:
$$
y_j^{(n_j)}(x_j)+\di\sum_{\mu=0}^{n_j-2} \Big(\frac{\nu_{\mu
j}}{x_j^{n_j-\mu}}+ q_{\mu j}(x_j)\Big)y_j^{(\mu)}(x_j) =\la
y_j(x_j),\quad x_j\in (0, l_j),\quad j=\overline{1,p},               \eqno(1)
$$
where $\la$ is the spectral parameter, $q_{\mu j}(x_j)$ are
complex-valued integrable functions. We call $q_j=\{q_{\mu
j}\}_{\mu=\overline{0,n_j-2}}$ the potential on the edge $e_j$,
and we call $q=\{q_{j}\}_{j=\overline{1,p}}$ the potential on
the graph $T.$ 
In this paper we study inverse spectral problems for system (1).
We provide a procedure for constructing the solution of the inverse 
problem and prove its uniqueness. 

\medskip
{\bf 2. } Let us construct special fundamental systems of solutions
for higher order differential operators with regular singularities.
Consider the differential equation
$$
\ell y:= y^{(n)}+\di\sum_{j=0}^{n-2} \Big(\di\frac{\nu_j}{x^{n-j}}
+q_j(x)\Big)y^{(j)}=\lambda y, \quad x>0                           \eqno(2)
$$
on the half-line. Let $\mu_1,\ldots,\mu_n$ be the roots of the
characteristic polynomial
$$
\Delta(\mu)=\di\sum_{j=0}^{n} \nu_j \di\prod_{k=0}^{j-1} (\mu-k),
\quad \nu_n=1,\; \nu_{n-1}=0.
$$
It is clear that $\mu_1+\ldots+\mu_n=n(n-2)/2.$ For definiteness, we assume
that $\mu_k-\mu_j\ne sn$ $(s=0,\pm 1,\pm 2,\ldots)$; $Re\,\mu_1
< \ldots<Re\,\mu_n$, $\mu_k \ne 0,1,2,\ldots,n-3$ (the other cases require
minor modifications). Let $\theta_n=n-1-Re(\mu_n-\mu_1)$, $q_{0j}(x)=
q_j(x)$ for $x\ge 1,$ and $q_{0j}(x)=q_j(x)x^{min(\theta_n-j,0)}$ for
$x\le 1$ and assume that $q_{0j}(x)\in L(0,\infty),$ $j=\overline{0,n-2}.$

First of all, we consider the following differential equation 
without spectral parameter:
$$
\ell_0 y:=y^{(n)}+\di\sum_{j=0}^{n-2}
\di\frac{\nu_j}{x^{n-j}}\,y^{(j)}=y.                              \eqno(3)
$$
Let $x=r\exp(i\varphi),\; r>0,\; \varphi\in(-\pi,\pi],\; x^\mu=
\exp(\mu(\ln r+i\varphi))$ and $\Pi_{-}$ be the $x$-plane with the cut
along the semi-axis $x\le 0.$ Take numbers $c_{j0},$ $j=\overline{1,n}$
from the condition
$$
\di\prod_{j=1}^{n} c_{j0}=(\det[\mu_j^{\nu-1}]_{j,\nu=\overline{1,n}})^{-1}.
$$
Then the functions
$$
C_j(x)=x^{\mu_j} \di\sum_{k=0}^{\infty} c_{jk}x^{nk},\quad c_{jk}(x)
=c_{j0}\Big(\di\prod_{s=1}^{k} \Delta(\mu_j+sn)\Big)^{-1}         \eqno(4)
$$
are solutions of (3), and $\det[C_j^{(\nu-1)}(x)]_{j,\nu=\overline{1,n}}
\equiv 1.$ Moreover, the functions $C_j(x)$ are analytic in $\Pi_{-}$.
Denote
$$
\ee_k=\exp\Big(\di\frac{2\pi i(k-1)}{n}\Big),\quad S_\nu=\Big\{x:\arg x
\in\Big(\di\frac{\nu\pi}{n},\di\frac{(\nu+1)\pi}{n}\Big)\Big\},
$$
$$
S_1^{*}=\bar{S}_{n-1}, \quad S_k^{*}=\bar{S}_{n-2k+1}\cup\bar{S}_{n-2k+2},
\; k=\overline{2,n};
$$
$$
Q_k=\Big\{x:\arg x\in\Big[\max\Big(-\pi,(-2k+2)\di\frac{\pi}{n}\Big),
\min\Big(\pi,(2n-2k+2)\di\frac{\pi}{n}\Big)\Big]\Big\},\quad
k=\overline{1,n}.
$$
For $x\in S_k^*$ equation (3) has the solutions $e_k(x),$ $k=
\overline{1,n}$ of the form 
$$
e_k^{(\nu-1)}(x)=\ee_k^{\nu}\exp(\ee_k x)z_{k\nu}(x),\; 
\nu=\overline{0,n-1},
$$ 
where $z_{k\nu}(x)$ are solutions of the integral equations
$$
z_{k\nu}(x)=1+\di\frac{1}{n} \di\int_x^{\infty} \Big(\di\sum_{j=1}^{n}
\ee_j^{\nu+1}\ee_k^{-\nu}\exp((\ee_k -\ee_j)(t-x))\Big)
\Big(\di\sum_{m=0}^{n-2} \nu_m\ee_k^{m}t^{m-n}z_{km}(t)\Big)\,dt
$$
(here $\arg t=\arg x,\; |t|>|x|$). Using the fundamental system of 
solutions $\{C_j(x)\}_{j=\overline{1,n}}$ one can write
$$
e_k(x)=\di\sum_{j=1}^{n} \beta_{kj}^0 C_j(x).                      \eqno(5)
$$
In particular this gives the analytic continuation for $e_k(x)$ on 
$\Pi_{-}.$

\smallskip
{\bf Lemma 1. }{\it The functions $\{e_k(x)\}_{k=\overline{1,n}},\;
x\in\Pi_{-}$ form a fundamental system of solutions of equation (3), and
$$
\det[e_k^{(\nu-1)}(x)]_{k,\nu=\overline{1,n}}=
\det[\ee_k^{\nu-1}]_{k,\nu=\overline{1,n}}.
$$
The asymptotics
$$
e_k^{(\nu-1)}(x)=\ee_k^{\nu-1}\exp(\ee_k x)(1+O(x^{-1})),
\quad |x|\to\infty,\; x \in Q_k.                                \eqno(6)
$$
is valid.}

\smallskip
We observe that the asymptotics (6) holds in the sectors $Q_k$ 
which are wider that the sectors $S_k^*$. Next we obtain 
connections between the Stokes multipliers $\beta_{kj}^0.$

\smallskip
{\bf Lemma 2. }{\it The following relations hold}
$$
\beta_{kj}^0=
\beta_{1j}^0 \ee_{k}^{\mu_j},\quad j,k=\overline{1,n},        \eqno(7)
$$
$$
\di\prod_{j=1}^n \beta_{1j}^0 =
(\det[\ee_k^{\mu_j}]_{k,j=\overline{1,n}})^{-1}
\det[\ee_k^{j-1}]_{k,j=\overline{1,n}}\ne 0.                  \eqno(8)
$$

\medskip
Indeed, for $\arg x\in (-\pi,\pi-\frac{2\pi s}{n})$ we have, 
by virtue of (4)-(5),
$$
e_k(\ee^s x)=
\di\sum_{j=1}^{n} \beta_{kj}^0(\ee^s)^{\mu_j}C_j(x).          \eqno(9)
$$
It is easily seen from construction of the functions $e_k(x)$ that
$e_1(\ee^s x)=e_{s+1}(x).$ Substituting (5) in this equality 
and comparing the corresponding coefficients, we obtain
(7). After that (8) becomes obvious.
\B

\smallskip
Now we consider the differential equation
$$
\ell_0 y=\lambda y=\rho^n y, \quad x>0.                      \eqno(10)
$$
It is evident that if $y(x)$ is a solution of (3), then
$y(\rho x)$ satifies (10). Define $C_j(x,\lambda)$ by
$$
C_j(x,\lambda)=\rho^{-\mu_j}C_{j}(\rho x)=x^{\mu_j}
\di\sum_{k=0}^{\infty} c_{jk}(\rho x)^{nk}.
$$
The functions $C_j(x,\lambda)$ are entire in $\lambda$, and
$\det[C_j^{(\nu-1)}(x,\lambda)]_{j,\nu=\overline{1,n}}\equiv 1.$
From Lemmas 1 and 2 we get the following theorem.

\medskip
{\bf Theorem 1. }{\it In each sector $S_{k_0}=\{\rho:\; \arg\rho
\in (\frac{k_0 \pi}{n},\frac{(k_0+1)\pi}{n})\}$ equation (10) has a
fundamental system of solutions $B_0=\{y_k(x,\rho)\}_{k=\overline{1,n}}$
such that $y_k(x,\rho)=y_k(\rho x),$
$$
|y_k^{(\nu)}(x,\rho)(\rho R_k)^{-\nu}
\exp(-\rho R_k x)-1| \le \di\frac{M_0}{|\rho|x},\;
\rho\in\bar{S}_{k_0},\; |\rho|x\ge 1,\; \nu=\overline{0,n-1},   \eqno(11)
$$
$$
\det[y_k^{(\nu-1)}(x,\rho)]_{k,\nu=\overline{1,n}}\equiv
\rho^{n(n-1)/2}\Omega,\quad \Omega:=
\det[R_k^{\nu-1}]_{k,\nu=\overline{1,n}}\ne 0,                  \eqno(12)
$$
$$
y_k(x,\rho)=\di\sum_{j=1}^{n} b_{kj}^0 \rho^{\mu_j}C_j(x,\lambda),
\quad b_{kj}^0=\beta_{j}^0 R_{k}^{\mu_j},\; \beta_{j}^0\ne 0,   \eqno(13)
$$
where the constant $M_0$ depens only on $\nu_j.$}

\medskip
The functions $y_k(x,\rho)$ are analogues of the Hankel functions
for the Bessel equation. Denote
$$
C_j^{*}(x,\lambda)=\det[C_k^{(\nu)}(x,\lambda)]_{\nu=\overline{0,n-2};
\,k=\overline{1,n}\setminus n-j+1},
$$
$$
y^*_j(x,\rho)=(-1)^{n-j}\Big(\rho^{(n-1)(n-2)/2}\Omega\Big)^{-1}
\det[y_k^{(\nu)}(x,\rho)]_{\nu=\overline{0,n-2};\,
k=\overline{1,n}\setminus j},
$$
$$
F_{k\nu}(\rho x)=\left \{ \begin{array}{ll}
R_{k}^{\nu}\exp(\rho R_k x),\; & |\rho|x>1,\\[3mm]
(\rho x)^{\mu_1-\nu},\; & |\rho|x\le 1,
\end{array} \right.
\quad F_{k}^*(\rho x)=\left \{ \begin{array}{ll}
\exp(-\rho R_k x),\; & |\rho|x>1,\\[3mm]
(\rho x)^{n-1-\mu_n},\; & |\rho|x\le 1,
\end{array} \right.
$$
$$
U_{k\nu}^{0}(x,\rho)=y_k^{(\nu)}(x,\rho)(\rho^{\nu}F_{k\nu}(\rho x))^{-1},
\quad U_{k}^{0,*}(x,\rho)=y^{*}_k(x,\rho)(F^*_{k}(\rho x))^{-1},
$$
$$
g(x,t,\lambda)=\di\sum_{j=1}^{n} (-1)^{n-j}C_{j}(x,\lambda)
C_{n-j+1}^*(t,\lambda)=\di\frac{1}{\rho^{n-1}}\di\sum_{j=1}^{n}
y_{j}(x,\rho) y_{j}^*(t,\rho).
$$
The function $g(x,t,\lambda)$ is the Green's function of the 
Cauchy problem $\ell_0 y-\lambda y=f(x),$ $y^{(\nu)}(0)=0,$ 
$\nu=\overline{0,n-1}.$ Using (11)-(13), we obtain
$$
|U_{k \nu}^{0}(x,\rho)|\le M_1,\quad |U_{k}^{0,*}(x,\rho)|
\le M_1,\quad x\ge 0,\quad \rho\in\bar{S}_{k_0},                     \eqno(14)
$$
$$
\left. \begin{array}{c}
|C_{j}^{(\nu)}(x,\lambda|\le M_2 |x^{\mu_j-\nu}|,\\[4mm]
\Big|\di\frac{\partial^{\nu}}{\partial x^{\nu}}g(x,t,\lambda)\Big|
\le M_2 \di\sum_{j=1}^{n} |x^{\mu_j-\nu} t^{n-1-\mu_j}|,\quad
|\rho x|\le C_0,\; t\le x,
\end{array} \right \}                                               \eqno(15)
$$
where $M_1$ depens on $\nu_j,$ and $M_2$ on $\nu_j$ and $C_0.$

\smallskip
Now we are going to construct fundamental systems of solutions
of equation (2). Denote
$$
J(\rho)=\di\sum_{m=0}^{n-2} J_m(\rho),
$$
$$
J_m(\rho)=|\rho|^{Re(\mu_1-\mu_n)}\di\int_0^{|\rho|^{-1}} t^{\theta_n-m}
|q_m(t)|\,dt +|\rho|^{m-n+1} \di\int_{|\rho|^{-1}}^{\infty} |q_m(t)|\,dt.
$$

{\bf Lemma 3. }{\it The following estimate holds}
$$
J(\rho)\le \di\frac{Q}{|\rho|},\; |\rho|\ge1, \quad
Q:= \di\sum_{m=0}^{n-2} \di\int_{0}^{\infty} |q_{0m}(t)|\,dt.
$$

Indeed, if $\theta_m-m\le 0$, then $Re(\mu_n-\mu_1)\ge n-m-1,$
and consequently
$$
J_m(\rho)\le |\rho|^{m-n+1}\Big(\di\int_0^{|\rho|^{-1}} t^{\theta_n-m}
|q_m(t)|\,dt + \di\int_{|\rho|^{-1}}^{\infty} |q_m(t)|\,dt\Big)
\le |\rho|^{m-n+1} \di\int_0^{\infty} |q_{0m}(t)|\,dt.
$$
If $\theta_m-m>0,$ then
$$
J_m(\rho)\le |\rho|^{m-n+1} \di\int_0^{\infty} |q_{m}(t)|\,dt
\le |\rho|^{m-n+1} \di\int_0^{\infty} |q_{0m}(t)|\,dt.
$$
Hence $J(\rho)\le Q|\rho|^{-1},\;|\rho|\ge 1,$ and Lemma 3 is proved.

We now construct the functions $S_j(x,\lambda),$ $j=\overline{1,n}$
from the system of integral equations
$$
S_j^{(\nu)}(x,\lambda)=C_j^{(\nu)}(x,\lambda)-\di\int_0^x
\di\frac{\partial^{\nu}}{\partial x^{\nu}}\,g(x,t,\lambda)
\Big(\di\sum_{m=0}^{n-2} q_m(t)S_j^{(m)}(t,\lambda)\Big)\,dt,
\;\nu=\overline{0,n-1}.                                          \eqno(16)
$$
By (15), system (16) has a unique solution; moreover the functions
$S_j^{(\nu)}(x,\lambda)$ are entire in $\lambda$ for each $x>0,$ the
functions $\{S_j(x,\lambda)\}_{j=\overline{1,n}}$ form a fundamental
system of solutions for equation (2),
$\det[S_j^{(\nu-1)}(x,\lambda)\}_{j,\nu=\overline{1,n}}\equiv 1,$ and
$$
S_j^{(\nu)}(x,\lambda)=O(x^{\mu_j-\nu}),\quad (S_j(x,\lambda)-
C_j(x,\lambda))x^{-\mu_j}=o(x^{\mu_n-\mu_1}),\quad x\to 0.      \eqno(17)
$$
Let $S_{k_0,\alpha}=\{\rho:\; \rho\in S_{k_0},\; |\rho|>\alpha\},
\; \rho_0=2M_1 Q +1.$ For $k=\overline{1,n}$, 
$\rho\in\bar{S}_{k_0,\rho_0}$ we consider the system of integral 
equations
$$
U_{k\nu}(x,\rho)=U^0_{k\nu}(x,\rho)+ \di\sum_{m=0}^{n-2}
\di\int_0^\infty A_{k\nu m}(x,t,\rho)U_{km}(t,\rho)\,dt,
\; x\ge0,\; \nu=\overline{0,n-1},                               \eqno(18)
$$
where
$$
A_{k\nu m}(x,t,\rho)=\di\frac{q_m(t)F_{km}(\rho t)}{\rho^{n-1-m}
F_{k\nu}(\rho x)}
\left \{ \begin{array}{ll}
-\di\sum_{j=1}^{k} F_{j\nu}(\rho x)U^0_{j\nu}(x,\rho)
F_j^*(\rho t)U^{0,*}_{j}(t,\rho), & t\le x,\\[5mm]
\di\sum_{j=k+1}^n F_{j\nu}(\rho x)U^0_{j\nu}(x,\rho)
F_j^*(\rho t) U^{0,*}_{j}(t,\rho), & t> x.
\end{array} \right.
$$
Using (14) and Lemma 3, we obtain
$$
\di\sum_{m=0}^{n-2} \di\int_0^\infty |A_{k\nu m}(x,t,\rho)|\,dt
\le M_1 J(\rho)\le \di\frac{M_1 Q}{|\rho|}.
$$
Consequently, system (18) with $\rho\in\bar{S}_{k_0,\rho_0}$
has a unique solution, and uniformly in $x\ge 0,$
$$
U_{k\nu}(x,\rho)-U^0_{k\nu}(x,\rho)=O(\rho^{-1}),
\quad \rho\in\bar{S}_{k_0,\rho_0}.                             \eqno(19)
$$

{\bf Theorem 2. }{\it For $x>0,\; \rho\in\bar{S}_{k_0,\rho_0}$
there exists an fundamental system of solutions of equation (2),
$B=\{Y_k(x,\rho)\}_{k=\overline{1,n}}$ of the form
$$
Y_k^{(\nu)}(x,\rho)=\rho^{\nu}F_{k\nu}(\rho x)U_{k\nu}(x,\rho),
$$
where the functions $U_{k\nu}(x,\rho)$ are solution of (18),
and (19) is true.

The function $Y_k^{(\nu)}(x,\rho)$ considered for each $x>0,$ are analytic
in $\rho\in {S}_{k_0,\rho_0},$ continuous in $\rho\in\bar{S}_{k_0,\rho_0}$
and $\det[Y_k^{(\nu-1)}(x,\rho)]_{k,\nu=\overline{1,n}}=\rho^{n(n-1)/2}
\Omega(1+O(\rho^{-1}))$ as $|\rho|\to\infty.$ The functions $Y_k(x,\rho)$
satisfy the equality
$$
Y_{k}(x,\rho)= y_{k}(x,\rho)-\di\frac{1}{\rho^{n-1}} \di\int_0^x
\Big(\di\sum_{j=1}^{k} y_{j}(x,\rho)y_{j}^{*}(t,\rho)\Big)\Big(
\di\sum_{m=0}^{n-2} q_{m}(t)Y_{k}^{(m)}(t,\rho)\Big)\,dt
$$
$$
+\di\frac{1}{\rho^{n-1}} \di\int_x^\infty \Big(\di\sum_{j=k+1}^{n}
y_{j}(x,\rho)y_{j}^{*}(t,\rho)\Big)\Big(\di\sum_{m=0}^{n-2}
q_{m}(t)Y_{k}^{(m)}(t,\rho)\Big)\,dt.
$$
Moreover, one has the representation
$$
Y_{k}(x,\rho)=\di\sum_{j=1}^{n} b_{kj}(\rho)S_{j}(x,\lambda),     \eqno(20)
$$
where}
$$
b_{kj}(\rho)=b^0_{kj}(\rho)\rho^{\mu_j}(1+O(\rho^{-1})),
\quad |\rho|\to\infty,\; \rho\in\bar{S}_{k_0,\rho_0}.             \eqno(21)
$$

The only part of the theorem that needs proof is the asymptotic 
formula (21). Let $\rho$ be fixed, $x\le |\rho|^{-1}.$ Then (13) 
and (20) become
$$
\left.\begin{array}{c}
U_{k0}^{0}(x,\rho)=\di\sum_{j=1}^{n} b_{kj}^0(\rho x)^{\mu_j-\mu_1}
\hat{C}_{j}(x,\lambda), \\[5mm]
U_{k0}(x,\rho)=\di\sum_{j=1}^{n} b_{kj}(\rho)(\rho)^{-\mu_1}
x^{\mu_j-\mu_1} \hat{S}_{j}(x,\lambda),
\end{array}\right\}                                               \eqno(22)
$$
where
$$\hat{C}_{j}(x,\lambda)=x^{-\mu_j}C_{j}(x,\lambda),\;
\hat{S}_{j}(x,\lambda)=x^{-\mu_j}S_{j}(x,\lambda),\quad
\hat{S}_{j}(0,\lambda)=\hat{C}_{j}(0,\lambda)=c_{j0}\ne 0.
$$
It follows from (22) that
$$
U_{k0}(x,\rho)-U_{k0}^{0}(x,\rho)=\di\sum_{j=1}^{n} \Big(b_{kj}(\rho)
\rho^{-\mu_1}-b_{kj}^0 \rho^{\mu_j-\mu_1}\Big)x^{\mu_j-\mu_1}
\hat{S}_{j}(x,\lambda)
$$
$$
+\di\sum_{j=1}^{n} b_{kj}^{0}(\rho x)^{\mu_j-\mu_1}
(\hat{S}_{j}(x,\lambda)-\hat{C}_{j}(x,\lambda)).                 \eqno(23)
$$
Denote
$$
\left.\begin{array}{c}
{\cal F}_{k1}(x,\rho)=U_{k0}(x,\rho)-U_{k0}^{0}(x,\rho),\\[3mm]
{\cal F}_{k,s+1}(x,\rho)=\Big({\cal F}_{ks}(x,\rho)-{\cal F}_{ks}(0,\rho)
\hat{S}_{s}(x,\lambda)c_{s0}^{-1}\Big) x^{\mu_s-\mu_{s+1}},\;
s=\overline{1,n-1}.
\end{array} \right \}                                            \eqno(24)
$$

{\bf Lemma 4. }{\it The following relations hold}
$$
(b_{ks}(\rho)\rho^{-\mu_1}-b_{ks}^0\rho^{\mu_s-\mu_1})
c_{s0}={\cal F}_{ks}(x,\rho),\quad s=\overline{1,n},             \eqno(25)
$$
$$
{\cal F}_{ks}(x,\rho)=((U_{k0}(x,\rho)-U_{k0}^{0}(x,\rho))
-\di\sum_{j=1}^{s-1} (b_{kj}(\rho)\rho^{-\mu_1}-b_{kj}^0
\rho^{\mu_j-\mu_1}) x^{\mu_j-\mu_1}\hat{S}_{j}(x,\lambda))
x^{\mu_1-\mu_{s}},                                               \eqno(26)
$$
$$
s=\overline{1,n}.
$$

{\it Proof. } When $s=1$ equality (25) follows from (23)
for $x=0,$ while (26) is obviously true. Suppose now that
(25) and (26) have been proved for $s=1,\ldots,N-1.$ Then
$$
\Big((U_{k0}(x,\rho)-U_{k0}^{0}(x,\rho))-\di\sum_{j=1}^{N-1}
(b_{kj}(\rho)\rho^{-\mu_1}-b_{kj}^0\rho^{\mu_j-\mu_1})x^{\mu_j-\mu_1}
\hat{S}_{j}(x,\lambda))\Big)x^{\mu_1-\mu_{N}}
$$
$$
=\Big((U_{k0}(x,\rho)-U_{k0}^{0}(x,\rho))-\di\sum_{j=1}^{N-2}
(b_{kj}(\rho)\rho^{-\mu_1}-b_{kj}^0\rho^{\mu_j-\mu_1})x^{\mu_j-\mu_1}
\hat{S}_{j}(x,\lambda))\Big)x^{\mu_1-\mu_{N-1}} x^{\mu_{N-1}-\mu_N}
$$
$$
-(b_{k,N-1}(\rho)\rho^{-\mu_1}-b_{k,N-1}^0 \rho^{\mu_{N-1}-\mu_1})
\hat{S}_{N-1}(x,\lambda) x^{\mu_{N-1}-\mu_N}={\cal F}_{kN}(x,\rho),
$$
which gives (26) for $s=N.$ We now write (23) as
$$
{\cal F}_{kN}(x,\rho)=\di\sum_{j=N}^{n} (b_{kj}(\rho)\rho^{-\mu_1}
-b_{kj}^0 \rho^{\mu_j-\mu_1})x^{\mu_j-\mu_N}\hat{S}_{j}(x,\lambda)
$$
$$
+\di\sum_{j=1}^{n} (b_{kj}^0(\rho x)^{\mu_j-\mu_1}
(\hat{S}_{j}(x,\lambda)-\hat{C}_{j}(x,\lambda)) x^{\mu_1-\mu_N}.
$$
Hence, using (17), we infer ${\cal F}_{kN}(0,\rho)=(b_{kN}(\rho)
\rho^{-\mu_1}-b_{kN}^0\rho^{\mu_N-\mu_1})c_{N0},$ which gives (25)
for $s=N.$ Lemma 4 is proved.

\smallskip
Now we write (18) for $\nu=0$ as
$$
{\cal F}_{k1}(x,\rho)=\di\frac{1}{\rho^{n-1}}\big(-\di\int_{0}^x
(\di\sum_{j=1}^{n} (U_{j0}^0(x,\rho) U_{j}^{0,*}(t,\rho))
(\rho t)^{n-1-\mu_n}V_k(t,\rho)\,dt
$$
$$
+\di\int_{0}^\infty \Big(\di\sum_{j=k+1}^{n} U_{j0}^0(x,\rho)
U_{j}^{0,*}(t,\rho))F_j^*(\rho t)\Big)V_k(t,\rho)\,dt\Big),      \eqno(27)
$$
where
$$
V_k(t,\rho)=\di\sum_{m=0}^{n-2} 
q_m(t)\rho^m F_{km}(\rho t)U_{km}(t,\rho).
$$
Since for $t\le x\le |\rho|^{-1}$ we have
$$
\di\sum_{j=1}^{n} U_{j0}^0(x,\rho)U_{j}^{0,*}(t,\rho)=
\rho^{\mu_n-\mu_1} x^{-\mu_1}t^{1-n+\mu_n} g(x,t,\lambda),
$$
it follows by way of (15) that
$$
\Big|\di\sum_{j=1}^{n} U_{j0}^0(x,\rho)U_{j}^{0,*}(t,\rho)\Big|
\le M_3 |(\rho x)^{\mu_n-\mu_1}|,\quad 0\le t\le x\le |\rho|^{-1}.  \eqno(28)
$$

{\bf Lemma 5. }{\it The following relations hold}
$$
{\cal F}_{ks}(0,\rho)=\rho^{\mu_s-\mu_1-n+1}c_{s0}
\di\int_{0}^{\infty} \Big(\di\sum_{j=k+1}^{n} b_{js}^0
F_{j}^*(\rho t)U_{j}^{0,*}(t,\rho)\Big)V_k(t,\rho)\,dt,             \eqno(29)
$$
$$
{\cal F}_{ks}(x,\rho)=\di\frac{1}{\rho^{n-1}}\Big(-x^{\mu_1-\mu_s}
\di\int_{0}^x \Big(\di\sum_{j=1}^{n} U_{j0}^0(x,\rho) U_{j}^{0,*}(t,\rho)
\Big)(\rho t)^{n-1-\mu_n}V_k(t,\rho)\,dt
$$
$$
-\di\sum_{\ell=1}^{s-1} x^{\mu_{\ell}-\mu_s} \di\int_{0}^\infty 
\Big(\di\sum_{j=k+1}^{n} b_{j\ell}^{0}\rho^{\mu_{\ell}-\mu_1}
(\hat{S}_\ell(x,\lambda)-\hat{C}_\ell(x,\lambda))F_{j}^*(\rho t)
U_{j}^{0,*}(t,\rho)\Big)V_k(t,\rho)\,dt
$$
$$
+\di\int_{0}^\infty \Big(\di\sum_{j=k+1}^{n} \Big(\di\sum_{\xi=s}^{n}
b_{j \xi}^{0} \rho^{\mu_{\xi}-\mu_1}\hat{C}_\xi(x,\lambda)\Big)
F_{j}^*(\rho t) U_{j}^{0,*}(t,\rho)\Big)V_k(t,\rho)\,dt \Big),
\quad x\le |\rho|^{-1}.                                           \eqno(30)
$$

{\it Proof. } For $s=1,$ (29) and (30) follow from (27),
in view of (22). Suppose now that (29) and (30) have been
proved for $s=1,\ldots,N.$ Then, using (24), we calculate
$$
{\cal F}_{k,N+1}(x,\rho)=({\cal F}_{kN}(x,\rho)-{\cal F}_{kN}(0,\rho)
\hat{S}_N(x,\lambda)c_{N0}^{-1}) x^{\mu_N-\mu_{N+1}}
$$
$$
=\di\frac{1}{\rho^{n-1}}\Big(-x^{\mu_1-\mu_{N+1}} \di\int_{0}^x
\Big(\di\sum_{j=1}^{n} (U_{j0}^0(x,\rho) U_{j}^{0,*}(t,\rho)\Big)
(\rho t)^{n-1-\mu_n}V_k(t,\rho)\,dt
$$
$$
-\sum_{\ell=1}^{N-1} x^{\mu_{\ell}-\mu_{N+1}} \di\int_{0}^\infty
\Big(\di\sum_{j=k+1}^{n} b_{j\ell}^{0} \rho^{\mu_{\ell}-\mu_1}
(\hat{S}_\ell(x,\lambda)-\hat{C}_\ell(x,\lambda))
F_{j}^*(\rho t) U_{j}^{0,*}(t,\rho)\Big)V_k(t,\rho)\,dt
$$
$$
+x^{\mu_{N}-\mu_{N+1}} \di\int_{0}^\infty \Big(\di\sum_{j=k+1}^{n}\Big(
\di\sum_{\xi=N}^{n} b_{j\xi}^{0} \rho^{\mu_{\xi}-\mu_1}x^{\mu_{\xi}-\mu_N}
\hat{C}_\xi(x,\lambda)-b_{jN}^{0}\rho^{\mu_{N}-\mu_{N+1}}\hat{C}_N(x,\lambda)
$$
$$
-b_{jN}^{0} \rho^{\mu_{N}-\mu_{N+1}}(\hat{S}_N(x,\lambda)-\hat{C}_N(x,
\lambda))\Big)F_{j}^*(\rho t) U_{j}^{0,*}(t,\rho))V_k(t,\rho)\,dt\Big)
$$
$$
=\di\frac{1}{\rho^{n-1}}\Big(-x^{\mu_1-\mu_{N+1}} \di\int_{0}^x
\Big(\di\sum_{j=1}^{n} U_{j0}^0(x,\rho) U_{j}^{0,*}(t,\rho)\Big)
(\rho t)^{n-1-\mu_n}V_k(t,\rho)\,dt
$$
$$
-\di\sum_{\ell=1}^{N} x^{\mu_{\ell}-\mu_{N+1}} \di\int_{0}^\infty
\Big(\di\sum_{j=k+1}^{n} b_{j\ell}^{0} \rho^{\mu_{\ell}-\mu_1}
(\hat{S}_\ell(x,\lambda)-\hat{C}_\ell(x,\lambda))F_{j}^*(\rho t)
U_{j}^{0,*}(t,\rho)\Big)V_k(t,\rho)\,dt
$$
$$
+\di\int_{0}^\infty \Big(\di\sum_{j=k+1}^{n}\Big(\di\sum_{\xi=N+1}^{n}
b_{j \xi}^{0} \rho^{\mu_{\xi}-\mu_1} x^{\mu_{s}-\mu_{N+1}}
\hat{C}_\xi(x,\lambda)\Big)F_{j}^*(\rho t) U_{j}^{0,*}(t,\rho)
\Big)V_k(t,\rho)\,dt\Big),
$$
i.e. (30) is valid for $s=N+1.$ We now let $x\to 0$ in
(30) for $s=N+1$, using (28), to obtain (29) for
$s=N+1.$ This proves Lemma 3.

\smallskip
It follows from (25) and (29) that
$$
b_{ks}(\rho)\rho^{-\mu_s}-b_{ks}^0=\di\frac{1}{\rho^{n-1}}
\di\int_{0}^\infty \Big(\di\sum_{j=k+1}^{n} b_{js}^{0}
F_{j}^*(\rho t)U_{j}^{0,*}(t,\rho)\Big)V_k(t,\rho)\,dt.        \eqno(31)
$$
Using (31), (14), (19) and Lemma 3, we obtain
$$
b_{ks}(\rho)\rho^{-\mu_s}-b_{ks}^0=O(J(\rho))=O(\rho^{-1}),
\quad |\rho|\to\infty,\; \rho\in \bar S_{k_0},
$$
i.e. (21) is valid. Theorem 2 is proved.

\medskip
{\bf 3.} Using constructed fundamental systems of solutions
on each edge we can study spectral properties of systems on graphs
and solve the inverse spectral problem. 
Let $\{\xi_{kj}\}_{k=\overline{1,n_j}}$
be the roots of the characteristic polynomial
$$
\de_j(\xi)=\sum_{\mu=0}^{n_j} \nu_{\mu j} \prod_{k=0}^{\mu-1}
(\xi-k), \quad \nu_{n_j,j}:=1,\; \nu_{n_j-1,j}:=0.
$$
For definiteness, we assume that $\xi_{kj}-\xi_{mj}\ne sn_j,
s\in{\bf Z},$ $Re\,\xi_{1j}<\ldots <Re\,\xi_{n_j,j},$
$\xi_{kj}\ne\overline{0,n_j-3}$ (other cases require minor
modifications). We set $\theta_j:=n_j-1-Re\,(\xi_{n_j,j}-\xi_{1j}),$
and assume that the functions $q_{\mu j}^{(\nu)}(x_j),$
$\nu=\overline{0,\mu-1},$ are absolutely continuous, and
$q_{\mu j}^{(\mu)}(x_j)x_j^{\theta_j}\in L(0,l_j).$

Fix $j=\overline{1,p}.$ Let the numbers $c_{kj0},\; k=\overline{1,n_j},$
be such that
$$
\prod_{k=1}^{n_j} c_{kj0}=
\Big(\det[\xi_{kj}^{\nu-1}]_{k,\nu=\overline{1,n_j}}\Big)^{-1}.
$$
According to results of the previous section one can construct 
the fundamental systems of solutions
$\{S_{kj}(x_j,\la)\}_{k=\overline{1,n_j}}$ of Eq. (1) on the edge
$e_j$ such that the functions $S_{kj}^{(\nu)}(x_j,\la),\;
\nu=\overline{0,n_j-1},$ are entire in $\la,$ and for each fixed
$\la,$ and $x_j\to 0,$ $S_{kj}(x_j,\la)\sim c_{kj0} x_j^{\xi_{kj}}$
Consider the linear forms
$$
U_{j\nu}(y_j)=\sum_{\mu=0}^{\nu}\ga_{j\nu\mu}y_j^{(\mu)}(l_j), \;
j=\overline{1,p},\; \nu=\overline{0,n_j-1},
$$
where $\ga_{j\nu\mu}$ are complex numbers,
$\ga_{j\nu}:=\ga_{j\nu\nu}\ne 0.$
Denote $\langle n\rangle:=(|n|+n)/2,$ i.e. $\langle n\rangle=n$
for $n\ge 0$, and $\langle n\rangle=0$ for $n\le 0.$ Fix
$s=\overline{1,p},\; k=\overline{1,n_s-1}.$ Let $\Psi_{sk}=
\{\psi_{skj}\}_{j=\overline{1,p}}$ be solutions of Eq. (1) on
the graph $T$ under the boundary conditions
$$
\psi_{sks}(x_s,\la)\sim c_{ks0}x_s^{\xi_{ks}},\quad x_s\to 0,     \eqno(32)
$$
$$
\psi_{skj}(x_j,\la)=O(x_j^{\xi_{\langle n_j-k-1\rangle+2,j}}),
\quad x_j\to 0,\quad j=\overline{1,p},\;j\ne s,                  \eqno(33)
$$
and the matching conditions at the vertex $v_0$:
$$
U_{1\nu}(\psi_{sk1})=U_{j\nu}(\psi_{skj}), \;\;
j=\overline{2,p},\; \nu=\overline{0,k-1},\;n_j>\nu+1,            \eqno(34)
$$
$$
\di\sum_{j=1,\,n_j>\nu}^{p} U_{j\nu}(\psi_{skj})=0,
\quad \nu=\overline{k,n_s-1}.                                    \eqno(35)
$$
Define additionally $\psi_{sns}(x_s,\la):=S_{ns}(x_s,\la).$
Using the solutions $\{S_{\mu j}(x_j,\la)\}$, one can write
$$
\psi_{skj}(x_j,\la)=\di\sum_{\mu=1}^{n_j} M_{skj\mu}(\la)
S_{\mu j}(x_j,\la),\quad j=\overline{1,p},\quad k=\overline{1,n_s-1}, \eqno(36)
$$
where the coefficients $M_{skj\mu}(\la)$ do not depend on $x_j.$
It follows from (36) and the boundary condition (32) for the
Weyl-type solutions that
$$
\psi_{sks}(x_s,\la)=S_{ks}(x_s,\la)+
\sum_{\mu=k+1}^{n_s}M_{sk\mu}(\la)S_{\mu s}(x_s,\la),
\quad M_{sk\mu}(\la):=M_{sks\mu}(\la).                            \eqno(37)
$$
Denote
$$
M_{s}(\la)=[M_{sk\mu}(\la)]_{k,\mu=\overline{1,n_s}},\quad
M_{sk\mu}(\la):=\de_{k\mu}\quad \mbox{for}\quad k\ge \nu.
$$
The matrix $M_s(\la)$ is called the Weyl-type matrix with respect
to the boundary vertex $v_s$. The inverse problem is formulated as
follows. Fix $w=\overline{2,p}.$

\smallskip
{\bf Inverse problem 1.} Given $\{M_{s}(\la)\},\;
s=\overline{1,p}\setminus w$, construct $q$ on $T.$

\smallskip
Fix $s=\overline{1,p},\; k=\overline{1,n_s-1}.$ Substituting (36)
into boundary and matching conditions (32)-(35), we obtain a linear
algebraic system with respect to $M_{skj\mu}(\la).$ Solving this
system one gets
$M_{skj\mu}(\la)=\Delta_{skj\mu}(\la)/\Delta_{sk}(\la),$ where the
functions $\Delta_{skj\mu}(\la)$ and $\Delta_{sk}(\la)$ are entire
in $\la.$ In particular,
$$
M_{sk\mu}(\la)= \Delta_{sk\mu}(\la)/\Delta_{sk}(\la),\quad k\le\mu,
$$
where $\Delta_{sk\mu}(\la):=\Delta_{sks\mu}(\la).$

\medskip
{\bf 4. } In this section we obtain
a constructive procedure for the solution of Inverse problem 1 and
prove its uniqueness. First we consider auxiliary inverse problems
of recovering the differential operator on each each fixed edge.
Fix $s=\overline{1,p},$ and consider the following inverse problem
on the edge $e_s$.

\smallskip
{\bf IP(s). } Given the matrix $M_s$,
construct the potential $q_{s}$ on the edge $e_s$.

\smallskip
{\bf Theorem 3. }{\it Fix $s=\overline{1,p}.$ The specification of the
Weyl-type matrix $M_s$ uniquely determines the potential $q_s$ on the
edge $e_s$.}

\smallskip
We omit the proof since it is similar to that in [1, Ch.2].
Moreover, using the method of spectral mappings, one can get a
constructive procedure for the solution of the inverse problem
$IP(s)$. It can be obtained by the same arguments as for $n$-th
order differential operators on a finite interval (see [1, Ch.2]
for details).

\smallskip
Fix $j=\overline{1,p}.$ Let $\vv_{jk}(x_j,\la),$ $k=\overline{1,n_j},$
be solutions of equation (1) on the edge $e_j$ under the conditions
$$
\vv_{kj}^{(\nu-1)}(l_j,\la)=\de_{k\nu},\;\nu=\overline{1,k},
\qquad \vv_{kj}(x_j,\la)=O(x_j^{\xi_{n_j-k+1,j}}),\;x_j\to 0.
$$
We introduce the matrix
$m_j(\la)=[m_{jk\nu}(\la)]_{k,\nu=\overline{1,n_j}},$
where $m_{jk\nu}(\la):=\vv^{(\nu-1)}_{jk}(l_j,\la).$
The matrix $m_j(\la)$ is called the Weyl-type matrix with
respect to the internal vertex $v_0$ and the edge $e_j$.

\smallskip
{\bf IP[j]. } Given the matrix $m_j$, construct $q_{j}$ on the edge $e_j$.

This inverse problem is the classical one, since it is the inverse
problem of recovering a higher-order differential equation on a
finite interval from its Weyl-type matrix. This inverse problem
has been solved in [1], where the uniqueness theorem for this
inverse problem is proved. Moreover, in [1] an algorithm for the
solution of the inverse problem $IP[j]$ is given, and necessary
and sufficient conditions for the solvability of this inverse
problem are provided.

\smallskip
Fix $j=\overline{1,p}.$ Then for each fixed $s=\overline{1,p}\setminus j,$
$$
m_{j1\nu}(\la)=\di\frac{\psi_{s1j}^{(\nu-1)}(l_j,\la)}{\psi_{s1j}(l_j,\la)},
\quad \nu=\overline{2,n_j},                                                \eqno(38)
$$
$$
m_{jk\nu}(\la)=\di\frac{\det[\psi_{s\mu j}(l_j,\la),\ldots,
\psi_{s\mu j}^{(k-2)}(l_j,\la),
\psi_{s\mu j}^{(\nu-1)}(l_j,\la)]_{\mu=\overline{1,k}}}
{\det[\psi_{s\mu j}^{(\xi-1)}(l_j,\la)]_{\xi,\mu=\overline{1,k}}}\,,
\;2\le k<\nu\le n_j.                                                      \eqno(39)
$$

Now we are going to obtain a constructive procedure for the
solution of Inverse problem 1. 
For this purpose it is convenient to divide differential
equations into $m$ groups with equal orders. More precisely, let
$\om_1>\om_2>\ldots >\om_m>\om_{m+1}=1,$ $n_{p_{j-1}+1}=\ldots=n_{p_j}
:=\om_j,$ $j=\overline{1,m},$ $0=p_0<p_1<\ldots <p_m=p.$
Take $N$ such that $p_N=w.$

Suppose that we already found the potentials
$q_{s}$, $s=\overline{1,p}\setminus p_N$, on the edges $e_s$,
$s=\overline{1,p}\setminus p_N$. Then we calculate the functions
$S_{kj}(x_j,\la),$ $j=\overline{1,p}\setminus p_N;$
here $k=\overline{1,\om_i}$ for $j=\overline{p_{i-1}+1,p_i}.$

Fix $s=\overline{1,p_1}$ (if $N>1$), and $s=\overline{1,p_1-1}$ (if $N=1$).
Our goal now is to construct the Weyl-type matrix $m_{p_N}(\la).$
According to (38)-(39), in order to construct $m_{p_N}(\la)$ we have
to calculate the functions
$$
\psi_{skp_N}^{(\nu)}(l_{p_N},\la),
\quad k=\overline{1,\om_N-1},\; \nu=\overline{0,\om_N-1}.                 \eqno(40)
$$
We will find the functions (40) by the following steps.

1) Using (37) we construct the functions
$$
\psi_{sks}^{(\nu)}(l_s,\la),
\;k=\overline{1,\om_N-1},\;\nu=\overline{0,\om_1-1},                     \eqno(41)
$$
by the formula
$$
\psi_{sks}^{(\nu)}(l_s,\la)=S_{ks}^{(\nu)}(l_s,\la)+
\di\sum_{\mu=k+1}^{\om_1} M_{sk\mu}(\la)S_{\mu s}^{(\nu)}(l_s,\la).       \eqno(42)
$$

2) Consider a part of the matching conditions (34) on $\Psi_{sk}$. More
precisely, let $\xi=\overline{N,m},\;k=\overline{\om_{\xi+1},\om_\xi-1},\;
l=\overline{\xi,m},\;j=\overline{1,p_l-1}.$ Then, in particular, (34) yields
$$
U_{p_l,\nu}(\psi_{skp_l})=U_{j\nu}(\psi_{skj}),
\quad \nu=\overline{\om_{l+1}-1,\min(k-1,\om_l-2)}.                      \eqno(43)
$$
Using (43) we can calculate the functions
$$
\psi_{skj}^{(\nu)}(l_j,\la),\;\xi=\overline{N,m},\;
k=\overline{\om_{\xi+1},\om_\xi-1},\; l=\overline{\xi,m},\;
j=\overline{1,p_l},\; \nu=\overline{\om_{l+1}-1,\min(k-1,\om_l-2)}.      \eqno(44)
$$

3) It follows from (36) and the boundary conditions on $\Psi_{sk}$ that
$$
\psi_{skj}^{(\nu)}(l_j,\la)=\di\sum_{\mu=\max(\om_l-k+1)}^{\om_l}
M_{skj\mu}(\la)C_{\mu j}^{(\nu)}(l_j,\la),                                \eqno(45)
$$
$$
k=\overline{1,\om_1-1},\;l=\overline{1,m},\;
j=\overline{p_{l-1}+1,p_l}\setminus s,\;\nu=\overline{0,\om_l-1}.
$$
We consider only a part of relations (45). More precisely, let
$\xi=\overline{N,m},\; k=\overline{\om_{\xi+1},\om_\xi-1},\;
l=\overline{1,m},\;j=\overline{p_{l-1}+1,p_l},\; j\ne p_N,\;
j\ne s,\; \nu=\overline{0,\min(k-1,\om_l-2)}.$ Then
$$
\di\sum_{\mu=\max(\om_l-k+1)}^{\om_l}M_{skj\mu}(\la)C_{\mu j}^{(\nu)}(l_j,\la)
=\psi_{skj}^{(\nu)}(l_j,\la),\quad \nu=\overline{0,\min(k-1,\om_l-2)}.    \eqno(46)
$$
For this choice of parameters, the right-hand side in (46) are known,
since the functions (44) are known. Relations (46) form a linear algebraic
system $\sigma_{skj}$ with respect to the coefficients $M_{skj\mu}(\la).$
Solving the system $\sigma_{skj}$ we find the functions $M_{skj\mu}(\la).$  
Substituting them into (45), we calculate the functions
$$
\psi_{skj}^{(\nu)}(l_j,\la),\quad k=\overline{1,\om_N-1},\; l=\overline{1,m},
\;j=\overline{p_{l-1}+1,p_l}\setminus p_N,\; \nu=\overline{0,\om_l-1}.    \eqno(47)
$$
Note that for $j=s$ these functions were found earlier.

4) Let us now use the generalized Kirchhoff's conditions (35) for $\Psi_{sk}$.
Since the functions (47) are known, one can construct by (35) the functions
(40) for $k=\overline{1,\om_N-1},\;\nu=\overline{k,\om_N-1}.$ Thus, the
functions (40) are known for $k=\overline{1,\om_N-1},\; \nu=\overline{0,\om_N-1}.$

\smallskip
Since the functions (40) are known, we construct the Weyl-type matrix
$m_{p_N}(\la)$ via (38)-(39) for $j=p_N.$ Thus, we have obtained the solution of
Inverse problem 1 and proved its uniqueness, i.e. the following assertion holds.

\smallskip
{\bf Theorem 2. }{\it The specification of the Weyl-type matrices $M_s(\la),$
$s=\overline{1,p}\setminus p_N$, uniquely determines the potential $q$ on $T.$
The solution of Inverse problem 1 can be obtained by the following algorithm.}

{\bf Algorithm 1. }{\it Given the Weyl-type matrices
$M_s(\la),$ $s=\overline{1,p}\setminus p_N$.

1) Find $q_{s}$, $s=\overline{1,p}\setminus p_N$, by solving the
inverse problem $IP(s)$ for each fixed $s=\overline{1,p}\setminus p_N$.

2) Calculate $C_{kj}^{(\nu)}(l_j,\la),\; j=\overline{1,p}\setminus p_N$; here
$k=\overline{1,\om_i},\;\nu=\overline{0,\om_i-1}$ for $j=\overline{p_{i-1}+1,p_i}.$

3) Fix $s=\overline{1,p_1}$ (if $N>1$), and $s=\overline{1,p_1-1}$ (if $N=1$).
All calculations below will be made for this fixed $s.$
Construct the functions (41) via (42).

4) Calculate the functions (44) using (43).

5) Find the functions $M_{skj\mu}(\la),$ by solving the linear algebraic
systems $\sigma_{skj}$.

6) Construct the functions (40) using (35).

7) Calculate the Weyl-type matrix $m_{p_N}(\la)$ via (38)-(39) for $j=p_N$.

8) Construct the potential $q_{p_N}$
by solving the inverse problem $IP[j]$ for $j=p_N$.}

We note that inverse spectral problems for Sturm-Liouville operators on spatial 
networks were studied in [2-7], and for higher order operators in [8].

\medskip
{\bf Acknowledgment.} This work was supported by Grant 1.1436.2014K of the Russian
Ministry of Education and Science and by Grant 13-01-00134 of Russian Foundation for
Basic Research.

\begin{center}
{\bf REFERENCES}
\end{center}
\begin{enumerate}
\item[{[1]}] V.A. Yurko. Method of Spectral Mappings in the Inverse Problem Theory,
     Inverse and Ill-posed Problems Series. VSP, Utrecht, 2002.
\item[{[2]}] Belishev M.I. Boundary spectral Inverse Problem on a class of graphs
     (trees) by the BC method. Inverse Problems 20 (2004), 647-672.
\item[{[3]}] Yurko V.A. Inverse spectral problems for Sturm-Liouville operators
     on graphs. Inverse Problems, 21, no.3 (2005), 1075-1086.
\item[{[4]}] Brown B.M.; Weikard R. A Borg-Levinson theorem for trees. Proc. R.
     Soc. Lond. Ser. A Math. Phys. Eng. Sci. 461 (2005), no.2062, 3231-3243.
\item[{[5]}] Yurko V.A. Inverse problems for Sturm-Liouville operators on
     bush-type graphs. Inverse Problems, 25, no.10 (2009), 105008, 14pp.
\item[{[6]}] Yurko V.A. An inverse problem for Sturm-Liouville operators
     on A-graphs.  Applied Math. Letters, 23, no.8 (2010), 875-879.
\item[{[7]}] Yurko V.A. Inverse spectral problems for differential operators
     on arbitrary compact graphs. Journal of Inverse and Ill-Posed Proplems,
     18, no.3 (2010), 245-261.
\item[{[8]}] Yurko V.A. Inverse spectral problems for arbitrary order differential
     operators on noncompact trees. Journal of Inverse and Ill-Posed Problems,
     20, no.1 (2012), 111-132.
\end{enumerate}

\begin{tabular}{ll}
Name:             &   Yurko, Vjacheslav  \\
Place of work:    &   Department of Mathematics, Saratov State University \\
{}                &   Astrakhanskaya 83, Saratov 410012, Russia \\
E-mail:           &   yurkova@info.sgu.ru
\end{tabular}

\end{document}